

\magnification=\magstep1
\everymath={\displaystyle}

\parskip=\medskipamount
%
\centerline {\bf The $q$-Harmonic Oscillator and the Al-Salam and
Carlitz
polynomials\/}
\bigskip
\centerline {\it Dedicated to the Memory of Professor Ya. A.
Smorodinski\u\i \/}
\bigskip
\centerline {R. Askey \plainfootnote{\dag}{Department of
Mathematics,
University of Wisconsin, Madison, WI 53706, USA}
and S. K. Suslov \plainfootnote{\ddag}{Russian Scientific Center
\lq \lq
Kurchatov Institute", Moscow 123182, Russia}}
\medskip
{\parindent = 12pt \narrower\narrower\narrower
\noindent{\bf Abstract}. One more model of a $q$-harmonic
oscillator based on
the $q$-orthogonal polynomials of Al-Salam and Carlitz is
discussed.
The explicit form of $q$-creation and $q$-annihilation
operators, $q$-coherent states and an analog of the Fourier
transformation
are established. A connection of the kernel of this transform with
a family of
self-dual biorthogonal rational functions is observed. \par}
\bigskip
\centerline{\bf Introduction\/}
\medskip
Recent development in quantum groups has led to the so-called
$q$-harmonic
oscillators ( see, for example, Refs. [1--7] ). Presently known
models of $q$-oscillators are closely related with $q$-orthogonal
polynomials.
The $q$-analogs of boson operators have been introduced explicitly
in Refs. [3], [5] and [7], where the corresponding wave functions
were
constructed in terms of the continuous $q$-Hermite polynomials of
Rogers [8,9],
in terms of  the Stieltjes--Wigert polynomials [10,11] and in terms
of
$q$-Charlier polynomials of Al-Salam and Carlitz [12],
respectively.
The model related to the Rogers--Szeg\"o polynomials [13] was
investigated
in [1,6]. Here we introduce the explicit realization of
$q$-creation and
$q$-annihilation operators with the aid of another family of the
Al-Salam and
Carlitz  polynomials [12] when eigenvalues of the corresponding
$q$-Hamiltonian are unbounded. An attempt to unify $q$-boson
operators is
also made.

With a great deal of regret we dedicate this paper to the memory of
Yacob A.
Smorodinski\u\i, who suggested ten years ago that the special
case $q=1$ of this work is interesting and admits a generalization.
\bigskip
\centerline{\bf 1. The Al-Salam and Carlitz Polynomials\/}
\medskip
The aim of this Letter is to show that the
$q$-orthogonal polynomials $U_n^{(a)}(x;q)$ studied by Al-Salam and
Carlitz
are closely connected with the $q$-harmonic oscillator. To
emphasize
these relations we use the notation $u_n^{\mu}(x;q)= \mu ^{-n}
q^{-n(n-1)/2}
U_n^{(-\mu)}(x;q)$ for the Al-Salam and Carlitz polynomials. In our
notation they can be defined by the {\it three-term recurrence
relation \/}
of the form
$$
\mu q^{n}u_{n+1}^{\mu}(x;q) + (1-q^n)\,u_{n-1}^{\mu}(x;q)=
\left( x-(1-\mu )q^{n} \right) u_n^{\mu}(x;q)\,, \eqno (1)
$$
$u_0^{\mu}(x;q)=1\,,\,u_1^{\mu}(x;q)=\mu^{-1}(x-1+\mu )\,.$ These
polynomials are {\it orthogonal \/}
$$\int_{-\mu}^{1}u_m^{\mu}(x;q)\,u_n^{\mu}(x;q)\,d\alpha (x)=
(1+\mu) q^{-n(n-1)/2}\,{(q;q)_n \over \mu^n }\,\delta _{mn} \eqno
(2)
$$
with respect to a positive measure $d\alpha(x)$, where $\alpha(x)$
is
a step function with jumps
$$
{ q^k \over (-q\mu;q)_{\infty} (q,-q/\mu;q)_k }
$$
at the points $x=q^k, k=0,1,\dots ,$ and jumps
$$
{ \mu q^k \over (-q/\mu;q)_{\infty} (q,-q\mu;q)_k }
$$
at the points $x=-\mu q^k, k=0,1,\dots $ ( see, for example,
[12,14,15] ).
Here the usual notations are
$$
\eqalignno{
&(a;q)_n=\prod_{k=0}^{n-1}(1-aq^k)\,,\cr
&(a,b;q)_n=(a;q)_n(b;q)_n\,,&(3)\cr
&(a;q)_{\infty}=\lim_{n\to \infty }(a;q)_n\,.\cr }
$$
The orthogonality relation (2) can also be written in terms of
the $q$-{\it integral\/} of Jackson,
$$\int_{-\mu}^{1}u_m^{\mu}(x;q)\,u_n^{\mu}(x;q)\,\tilde
\rho(x)\,d_q x=
(1-q) d_n^2\,\delta _{mn}\,, \eqno (4)
$$
where
$$
\tilde \rho(x)= { (qx,-\mu^{-1}qx;q)_{\infty} \over
(q,-\mu,-q/\mu;q)_{\infty} }\, ;\;\, \mu>0\,,\,0<q<1 \eqno (5)
$$
and
$$
d^2_n=q^{-n(n-1)/2}\, { (q;q)_n \over \mu^n }.\eqno (6)
$$
For the definition of the $q$-integral, see [15].
The \lq \lq weight function" $\rho(s)= \tilde \rho(x)$ in (5) is a
solution
of the {\it Pearson-type equation\/} $\Delta(\sigma \rho) = \rho
\tau
\nabla x_1 \;$ with\;$ x (s) = q^s,\; \sigma (s) =(1-q^s)(\mu
+q^s)\;$ and
$\sigma (s) + \tau (s) \nabla x_1 (s) = \mu$. The polynomials
$y_n(s)=
u_n^{\mu}(x;q)$ satisfy the { \it hypergeometric-type difference
equation\/} in self-adjoint form,
$$
{ \Delta \over \nabla x_1(s) } \left[ \sigma (s)\, \rho (s)\,
{ \nabla y_n(s) \over \nabla x(s) } \right] + \lambda_n \,\rho (s)
\,y_n(s)=0\,,
$$
where
$$
\lambda_n = q^{3/2}\, { q^{-n}-1 \over (1-q)^2 }\,.
$$
Here $\Delta f(s)=f(s+1)-f(s)=\nabla f(s+1)$ and $x_1(s)=x(s+1/2)$.
( For details, see [16--19]. )  The orthogonality property (2) or
(4) can be
proved by using standard Sturm--Liouville-type arguments ( cf.
[16--19]).

The explicit form of the polynomials $u_n^{\mu }(x;q)$ is
$$
\eqalignno{
u_n^{\mu }(x;q) &= \,_2\varphi_1 \left(
q^{-n},\,x^{-1};\,0\,;\,q\,,\,-\,
{q\over\mu}\,x \right)
&(7) \cr
&=(-\mu^{-1})^n \,_2\varphi_1(q^{-n},-\mu
x^{-1};\,0\,;\,q\,,qx)\,,\;x=q^s\,.
\cr}
$$
It means
$u_n^{\mu}(x;q)=\left(-\mu^{-1}\right)^nu_n^{1/\mu}(-\mu^{-1}x;q)$.
In the limit $q\to 1$ it easy to see from (1) or (7) that
$$
\lim _{q\to 1}u_n^{(1-q)\mu} (q^s;q) = \,_2F_0 (-n,-s;-;-1/\mu ) =
c_n^{\mu}(s)\,, \eqno (8)
$$
where $c_n^{\mu}(x)$ are the Charlier polynomials.
\bigskip
\centerline{\bf 2. Model of $q$-Harmonic Oscillator\/}
\medskip
The Al-Salam and Carlitz polynomials $u_n^{\mu }(x;q)$ allow us to
consider
an interesting model of a $q$-oscillator ( cf. [7] ). We can
introduce
a $q$-{\it version of the wave functions\/} of the harmonic
oscillator as
$$
\psi_n(s)=\tilde \psi_n(x)=d^{-1}_n \left( \tilde \rho(x)|x|
\right) ^{1/2}\,
u_n^{\mu }(x;q)\,,\;x=q^s\,,  \eqno (9)
$$
where $\tilde \rho (x)$ and $d^2_n$ are defined in (5) and (6),
respectively.
These $q$-wave functions satisfy the orthogonality relation
$$
(1-q)^{-1} \int^1_{-\mu}\tilde \psi_n(x)\,\tilde
\psi_m(x)\,|x|^{-1}\,d_qx
= \delta_{nm} \,,\eqno (10)
$$
which is equivalent to (2) and (4).

The $q$-{\it annihilation\/} $b$ and $q$-{\it creation\/} $b^+$
operators
have the following explicit form
$$\eqalignno{
b=&\;(1-q)^{- {1 \over 2}} \left[\,\mu ^{1 \over 2} q^{-s} -
\sqrt{(1-q^{s+1})(\mu q^{-1}+q^s)} \,q^{-s}\,
e^{\partial _s} \right]\,,\cr
&&(11) \cr
b^+=&\;(1-q)^{- {1 \over 2}} \left[\,\mu ^{1 \over 2} q^{-s} -
e^{-\partial _s}
\sqrt {(1-q^{s+1}) (\mu q^{-1}+q^s)}\,q^{-s} \right]\,,\cr
}$$
where $\partial _s \equiv {d \over ds},\; e^{\alpha \partial_s}
f(s)=f(s+
\alpha)$. These operators are adjoint, $(b^+ \psi, \chi)=
(\psi, b\chi)$, with respect to the scalar product (10).
They satisfy the $q$-{\it commutation rule\/}
$$
b\,b^+-q^{-1}b^+b=1\eqno {(12)}
$$
and act on the $q$-wave functions defined in (9) by
$$
b\,\psi_n=\tilde e^{1/2}_n\,\psi_{n-1},\; \;
b^+\psi_n=\tilde e^{1/2}_{n+1}\psi_{n+1}\,,\eqno (13)
$$
where
$$
\tilde e_n={1-q^{-n} \over 1-q^{-1}}\,.
$$

The $q$-{\it Hamiltonian\/} $H=b^+b\,$ acts on the wave functions
(9) as
$$
H\psi_n=\tilde e_n\psi_n\,\eqno(15)
$$
and has the following explicit form
$$\eqalignno{
H &= (1-q)^{-1} \bigl[\,\mu q^{-2s}+(1-q^s)( \mu +q^s )\,q^{1-2s}-
&(16) \cr
 & \mu^{1 \over 2}q^{-2s} \sqrt{(1-q^{s+1})(\mu
q^{-1}+q^s)}\,e^{\partial _s}-
\mu^{1 \over 2}q^{2-2s} \sqrt{(1-q^s)(\mu
q^{-1}+q^{s-1})}\,e^{-\partial _s}
\bigr]\,. \cr }
$$
By factorizing the Hamiltonian ( or the difference equation for the
Al-Salam
and Carlitz polynomials ) we arrive at the explicit form (11) for
the $q$-boson
operators. The equations (13) are equivalent to the following
difference-differentiation formulas
$$\eqalign{
&\mu q^{-s-1} \Delta u_n^{\mu}(x;q)=(1-q^{-n})\,u_{n-1}^{\mu}(x;q)
\,,\cr
&q^{-s} \nabla \,[\, \rho(s) \, u_n^{\mu}(x;q)\,]= \rho(s)
\,u_{n+1}^{\mu}(x;q)
\,,\cr}
$$
respectively. Therefore, the main properties of the Al-Salam and
Carlitz
polynomials admit a simple group-theoretical interpretation in
terms of
the $q$-Heisenberg--Weyl algebra (12).
The symmetric case $\mu =1$ in the above formulas corresponds to
the
{ \it discrete $q$-Hermite polynomials\/} $H_n(x;q)$ [12,15].
\bigskip
\centerline{\bf 3. The $q$-Coherent States\/}
\medskip
For the model of the $q$-oscillator under discussion, by analogy
with [7]
we can construct explicitly the $q$-{\it coherent states\/}
$\mid \alpha \rangle$ defined by
$$
\eqalignno{
b\mid\alpha\rangle =&\, \alpha \mid\alpha\rangle \,, &(17)\cr
\mid\alpha\rangle =&\,
f_{\alpha}\sum_{n=0}^{\infty}{\alpha^n\psi_n(s) \over
(\tilde e_n!)^{1/2}}\,,
\quad \langle\alpha\mid\alpha\rangle =1\,,\cr}
$$
where
$$
\eqalign{
\tilde  e_n!=&\, \tilde e_1 \tilde e_2 \dots \tilde e_n =
q^{-n(n-1)/2}\,{ (q;q)_n \over (1-q)^n }\,,\cr
f_{\alpha}=&\,
\left(-(1-q)\mid\alpha\mid^2;q\right)^{-1/2}_{\infty}\,.
\cr}
$$
By using (9) one can obtain
$$
\mid\alpha\rangle =f_{\alpha} \left( \rho |q^s| \right)^{1 \over 2}
\sum_{n=0}^\infty u_n^{\mu}(x;q)\,q^{n(n-1)/2}\,
{t^n \over (q;q)_n}\,,\;t=\alpha \mu^{1 \over 2}
(1-q)^{1 \over 2}\,. \eqno(18)
$$
With the aid of the Brenke-type generating function [12,14] for
the Al-Salam and Carlitz polynomials,
$$
\sum_{n=0}^{\infty}u_n^{\mu}(x;q)\,q^{n(n-1)/2}\,
{t^n \over (q;q)_n}=
{(-t\,,\,t/ \mu ;q)_\infty \over
(xt/ \mu ;q)_\infty }\,,\; \left|t\,{x \over \mu}\right|<1\,,
\eqno(19)
$$
we arrive at the following {\it explicit form \/} for the
$q$-coherent states
$$
\mid\alpha\rangle =f_{\alpha}\left(\rho |q^s| \right)^{1 \over 2}
{\left(-\alpha (1-q)^{1/2}\mu ^{1/2},\alpha (1-q)^{1/2}\mu^{-1/2};q
\right)
_\infty \over
\left(\alpha (1-q)^{1/2} \mu^{-1/2}q^s;q \right)_\infty }
\,, \eqno{(20)}
$$
where
$\rho(s)=(q^{1+s},-\mu^{-1}q^{1+s};q)_\infty/(q,-\mu,-q/\mu;q)_\i
nfty $.
These coherent states are not orthogonal
$$
\langle\alpha\mid\beta\rangle ={ \left(-(1-q)\alpha^*\beta\,;q
\right)_\infty \over
\left(-(1-q) \mid\alpha\mid ^2\,,-(1-q) \mid\beta
\mid ^2;q \right)_\infty ^{1/2} } \,,
$$
where $*$ denotes the complex conjugate.
\bigskip
\centerline{\bf 4. Analog of the Fourier Transformation\/}
\medskip
To define an analog of the {\it Fourier transform \/} we begin,
in the spirit of Wiener's approach to the classical Fourier
transform [20]
( see also [7,21,22] ), by deriving the kernel of the form
$$\eqalignno{
K_t(x,y)&=\sum^\infty _{n=0}\,t^n \tilde \psi_n(x)\,\tilde
\psi_n(y) &(21) \cr
&= \left(\tilde\rho(x) \tilde\rho(y) |xy| \right)^{1 \over 2}
\sum_{n=0}^\infty u_n^{\mu}(x;q)\,u_n^{\mu}(y;q)\,q^{n(n-1)/2}\,
{(\mu t)^n \over (q;q)_n}\,.\cr }
$$
The series can be summed with the aid of the bilinear generating
function
of Al-Salam and Carlitz [12]
$$\eqalignno{
\sum_{n=0}^\infty u_n^{\mu _1}(x;q)\,u_n^{\mu _2}(y;q)\,
q^{n(n-1)/2}\,{t^n \over (q;q)_n}=
{(-t\,,\,t/ \mu_1\,,\,t/\mu_2 ;q)_\infty \over
(tx/ \mu_1\,,\,ty/\mu_2 ;q)_\infty }& \cr
\cdot \,_3\varphi_2\left(
\matrix
x^{-1}\,, \,y^{-1}\,,\, -qt^{-1}\, \\
\noalign{\smallskip}
q\mu_1 (tx)^{-1},\, q\mu_2 (ty)^{-1} \\
\endmatrix
;\,q\,,\, q \right)&  &(22) \cr }
$$
( the $_3\varphi_2$-series is terminating and
$max \,(\, |t/\mu_1|\,,\,|t/\mu_2|\, )<1$ ). The answer is
$$\eqalignno{
K_t(x,y)&= \left( \tilde\rho(x) \tilde\rho(y) |xy| \right)^
{1 \over 2}
{(t\,,\,t\,,\,-\mu t\,;q)_\infty \over
(tx\,,ty\,;q)_\infty }
\cdot \, _3\varphi_2\left(
\matrix
x^{-1}\,, \,y^{-1}\,,\,-q(\mu t)^{-1}\, \\
\noalign{\smallskip}
q(tx)^{-1}\,,\, q(ty)^{-1} \\
\endmatrix
;\, q\,,\,q \right) \,; \cr
&&(23)\cr
&= \left( \tilde\rho(x) \tilde\rho(y) |xy| \right)^{1 \over 2}
{(t\,,\,t\,,\,-\mu^{-1} t\,;q)_\infty \over
(-\mu^{-1}tx,-\mu^{-1}ty\,;q)_\infty }
\cdot \, _3\varphi_2\left(
\matrix
-\mu x^{-1}, -\mu y^{-1},-q\mu t^{-1}\, \\
\noalign{\smallskip}
-q\mu (tx)^{-1},-q\mu (ty)^{-1} \\
\endmatrix
;q,q \right) \cr }
$$
at $x=q^s$ and at $x=-\mu q^s$ for $s=0,1,\dots$, respectively.
\vfill
\eject

In view of (10) and (21),
$$
t^m \tilde \psi_m(x)=(1-q)^{-1}\int_{-\mu}^1 K_t(x,y)\,\tilde
\psi_m(y)\,
|y|^{-1}d_qy\,.\eqno(24)
$$
Letting $t=i$, we find that the $q$-wave functions (9) are
eigenfunctions
of the following \lq \lq $q$-Fourier transform'',
$$
i^m \tilde\psi_m(x)=(1-q)^{-1}\int_{-\mu}^1
K_i(x,y)\,\tilde\psi_m(y)\,
|y|^{-1}d_qy\,.\eqno(25)
$$
An easy corollary of (21) or (24) is
$$
(1-q)^{-1}\int_{-\mu}^1
K_t(x,y)\,K_{t'}(x',y)\,|y|^{-1}d_qy=K_{tt'}(x,x')\,.
\eqno(26)
$$
Putting $t=-t'=i$, we obtain the orthogonality relation of the
kernel,
$$
(1-q)^{-1}\int^1 _{-\mu} K_i(x,y)\,K_i^*(x',y)\,|y|^{-1}d_qy=
\delta_{xx'} \,,
\eqno (27)
$$
which implies the orthogonality of the rational functions (23) and
leads to
an inversion formula for the $q$-transformation (25).
In view of (8),
in the limit $q \to 1^-$ we get one of the \lq \lq discrete Fourier
transforms''
considered in [21].
\bigskip
\centerline{\bf 5. Some Biorthogonal Rational Functions\/}
\medskip
The rational functions (23) have appeared as the kernel of the
discrete
$q$-Fourier transform (25). They admit the following extension.
With the aid of the bilinear generating function (22) and
the orthogonality property of a special case of the $q$-Meixner
polynomials,
which are dual to the polynomials (7), we obtain the {\it
biorthogonality
relation\/},
$$
\int^1_{-\mu_2} u(x,y)\,v(x',y)\,
\tilde\rho(y)\,d_qy=(1-q)\,d^2_x\,\delta_{xx'}
\,,\eqno (28)
$$
for the $_3 \varphi _2 $-{\it rational functions\/} of the form
$$\eqalignno{
u(x,y)&= \,_3\varphi_2\left(
\matrix
x^{-1}\,,\, y^{-1}\,, \,-qt_1^{-1} \, \\
\noalign{\smallskip}
q\mu_1 (t_1x)^{-1}\,,\,q\mu_2 (t_1y)^{-1} \\
\endmatrix
;\, q,\,\,q \right) \,;\cr
&&(29)\cr
&= \,_3\varphi_2\left(
\matrix
-\mu_1x^{-1}\,,\, -\mu_2y^{-1}\,, \,-qt_2 \, \\
\noalign{\smallskip}
-qt_2x^{-1}\,,\,-qt_2y^{-1} \\
\endmatrix
;\, q,\,\,q \right) \cr }
$$
at $x=q^s$ and at $x=-\mu_1 q^s$ for $s=0,1,\dots$, respectively,
and
$$
v(x,y)=u(x,y) \vert _{ t_1 \leftrightarrow t_2 }\,; \qquad
t_1 t_2= \mu_1 \mu_2\,. \eqno (30)
$$
\vfill
\eject
Here,
$$\eqalign{
\tilde\rho(y)= &\; {\left(qy\,,\,-\mu_2^{-1}qy\,;\,q
\right)_{\infty} \over
\left(t_1 \mu_2^{-1}y \,,\,t_2 \mu_2^{-1}y\,;q
\right)_{\infty}}\,;\cr
= &\; {\left(qy\,,\,-\mu_2^{-1}qy\,;\,q \right)_{\infty} \over
\left(-t_1^{-1}y \,,\,-t_2^{-1}y\,;q \right)_{\infty}} \cr }
$$
at $x,x'=\{ q^s;\,s=0,1,\dots \} $ and at
$x,x'=\{ -\mu_1q^s;\,s=0,1,\dots \} $, respectively;
$$\eqalign{
\tilde\rho(y)= &\; {\left(qy\,,\,-\mu_2^{-1}qy\,;\,q
\right)_{\infty} \over
\left(t_1 \mu_2^{-1}y \,,\,-t_1^{-1}y\,;q \right)_{\infty}}\,;\cr
= &\; {\left(qy\,,\,-\mu_2^{-1}qy\,;\,q \right)_{\infty} \over
\left(-t_2^{-1}y \,,\,t_2\mu_2^{-1}y\,;q \right)_{\infty}} \cr }
$$
at $x,x'=\{ q^s,-\mu_1q^s \} $ and vice versa, respectively. The
squared norm
is
$$\eqalign{
d^2_x= &\; {\left( q, \,q,-\mu_1,-\mu_2,-q\mu_1^{-1},-q\mu_2^{-1};q
\right)_
\infty
\over
\left( -t_1,-t_2,\,t_1\mu_1^{-1},\,t_2\mu_1^{-1},\,t_1\mu_2^{-1},\,
t_2\mu_2^{-1};q \right)_\infty } \cdot
{\left( t_1\mu_1^{-1}x,\,t_2\mu_1^{-1}x;q \right)_{\infty} \over
\left( qx,-\mu_1^{-1}qx\,;q \right)_\infty } \, |x|^{-1}\,; \cr
= &\; {\left( q, \,q,-\mu_1,-\mu_2,-q\mu_1^{-1},-q\mu_2^{-1};q
\right)_\infty
\over
\left(
-t_1^{-1},-t_2^{-1},\,\mu_1t_1^{-1},\,\mu_1t_2^{-1},\,\mu_2t_1^{-
1},\,
\mu_2t_2^{-1};q \right)_\infty } \cdot
{\left( -t_1^{-1}x,\,-t_2^{-1}x;q \right)_{\infty} \over
\left( qx,-\mu_1^{-1}qx\,;q \right)_\infty } \, |x|^{-1} \cr }
$$
for $x=q^s$ and for $x=-\mu_1q^s$, respectively.

The functions (29)--(30) are self-dual and belong to  classical
biorthogonal
rational functions [23--27]. It is interesting to compare the
biorthogonality
relation (28) with
the orthogonality property for the big $q$-Jacobi polynomials [28],
which live at the same terminating $_3\varphi_2$-level.
\bigskip
\centerline{\bf 6. Concluding Remarks\/}
\medskip
In view of (11), it is natural to introduce operators
of the form
$$
a=\alpha(s)-\beta(s)\,e^{\partial}\,,\;
\;a^+=\alpha(s)-e^{-\partial}\beta(s)
$$
with two arbitrary functions $\alpha(s)$ and $\beta(s)$ and to
satisfy the
commutation rule $a\,a^+-qa^+a=1$. The result is
$$
\alpha(s+1)=q\alpha(s)\,,
$$
$$
(1-q)\alpha^2(s)+\beta^2(s)-q\beta^2(s-1)=1
$$
and we can choose $\alpha(s)=\varepsilon \,q^s$ and
$\beta^2(s)=\varepsilon
^2\,(q^{s+1}-\gamma)(q^s-\delta)$ with $(1-q)\gamma \delta
\varepsilon^2=1$.
Since
$$
(a^+\psi,\chi
)-(\psi,a\chi)=\sum_s\Delta[\beta(s-1)\,\psi^*(s-1)\,\chi(s)]\,,
$$
the corresponding operators are adjoint for the two different cases
considered
in [7] and in this Letter with $0<q<1$ and $q>1$, respectively.

For $\beta=constant$ we can try
$$
a=e^{\partial}\left( e^{\partial}-\alpha(s) \right) \,,\; \;
a^+=\left(e^{-\partial}-\alpha(s)\right)\,e^{-\partial}
$$
and obtain $aa^+-qa^+a=1-q$, when
$$
\alpha^2(s+1)=q\alpha^2(s)\,,\;\;\alpha(s+2)=q\alpha(s)\,,
$$
which is satisfied for \,$\alpha=\varepsilon \,q^{s/2}$. This case
has been
considered in [5].

Finally, the operators
$$
a=\varepsilon \, { e^{\gamma \partial}\alpha(s)+e^{-\gamma
\partial}\beta(s)
\over \alpha(s)-\beta(s) } \,,\; \;
a^+=\varepsilon \, { \alpha(s)e^{-\gamma \partial}+
\beta(s)e^{\gamma \partial}
\over \alpha(s)-\beta(s) }
$$
obey the $q$-commutation rule provided that $\alpha(s) \beta(s)
=\pm 1$ and
$\alpha(s+2\gamma)=q^{-1}\alpha(s)$. Therefore, $\alpha=q^{-s}$ for
$ \gamma =1/2$ and $\varepsilon^2=q^{1/2}(1-q)^{-1}$ ( cf. [3] ).

We can also introduce the operators
$$
a=\alpha^{-1}(s)- \varepsilon \, \beta^{-1}(s)\,e^{\partial}\,,\;
\;a^+=\alpha(s)- \varepsilon \, e^{\partial}\beta(s)
$$
and obtain $aa^+-qa^+a=1-q$ if
$$
\alpha(s+1)=q\alpha(s)\,,\;\;\beta(s+2)=q\beta(s)\,,
$$
so $\alpha =q^s$ and $\beta=q^{s/2}$. This leads to the
Rogers--Szeg\"o
polynomials [13] orthogonal on the unit circle ( see [1,6] ).
\bigskip
\centerline {\bf References }
\frenchspacing
\item{1.}  Macfarlane, A.J., {\it J. Phys. A: Math. Gen.\/} {\bf
22\/},
4581 (1989).

\item{2.}  Biedenharn, L.C.,  {\it J. Phys. A: Math. Gen.\/} {\bf
22\/},
L873 (1989).

\item{3.}  Atakishiyev, N.M. and Suslov, S.K., {\it Teor. i Matem.
Fiz.\/}
{\bf 85\/}, 64 (1990).

\item{4.}  Kulish, P.P. and Damaskinsky, E.V., {\it J. Phys. A:
Math. Gen.\/}
{\bf 23\/}, L415 (1990).

\item{5.}  Atakishiyev, N.M. and Suslov, S.K., {\it Teor. i Matem.
Fiz.\/}
{\bf 87\/}, 154 (1991).

\item{6.} Floreanini, R. and Vinet, L., {\it Lett. Math. Phys.\/}
{\bf 22\/}, 45 (1991).

\item{7.}  Askey, R. and Suslov, S.K., \lq The $q$-Harmonic
Oscillator and
an Analog of the Charlier Polynomials', Preprint No. 5613/1,
Kurchatov
Institute, Moscow 1993; {\it J. Phys. A: Math. Gen.\/}, submitted.

\item{8.} Rogers, L.J., {\it Proc. London Math. Soc.\/} {\bf 25\/},
318 (1894).

\item{9.} Askey, R. and Ismail, M.E.H., in {\it  Studies in Pure
Mathematics\/}
( P. Erd\"os, ed.), Birkh\"auser, Boston, Massachusetts, 1983, p.
55.

\item{10.}  Stieltjes, T.J., {\it Recherches sur les Fractions
Continues\/},
Annales de la Facult\'e des Sciences de Toulouse,
{\bf 8\/} (1894) 122 pp., {\bf 9\/} (1895), 47 pp.
Reprinted in {\it Oeuvres Compl\'etes\/}, vol. 2.

\item{11.}  Wigert S., {\it Arkiv f\"or Matematik, Astronomi och
Fysik\/},
{\bf 17\/}(18), 1 (1923).

\item{12.} Al-Salam, W.A. and Carlitz, L., {\it Math. Nachr.\/}
{\bf 30\/},
47 (1965).

\item{13.} Szeg\"o, G., {\it  Collected Papers \/}, Vol. 1
( R. Askey, ed.), Birkh\"auser, Basel, 1982, p. 795.

\item{14.} Chihara, T.S., {\it An Introduction to Orthogonal
Polynomials\/}, Gordon and Breach, New York, 1978.

\item{15.} Gasper, G. and Rahman, M., {\it Basic Hypergeometric
Series\/},
Cambridge Univ. Press, Cambridge, 1990.

\item{16.} Nikiforov, A.F. and Suslov, S.K., {\it Lett. Math.
Phys.\/}
{\bf 11\/}, 27 (1986).

\item{17.} Suslov, S.K., {\it Lett. Math. Phys.\/}  {\bf 14\/}, 77
(1987).

\item{18.} Suslov, S.K., {\it Russian Math. Surveys\/}, London
Math. Soc.
{\bf 44\/}, 227 (1989).

\item{19.} Nikiforov, A.F., Suslov, S.K., and Uvarov, V.B.,
{\it Classical Orthogonal Polynomials of a Discrete Variable\/},
Springer-Verlag, Berlin, Heidelberg, 1991.

\item{20.} Wiener, N., {\it The Fourier Integral and Certain of Its
Applications\/}, Cambridge University Press, Cambridge, 1933.

\item{21.} Askey, R., Atakishiyev, N.M., and Suslov, S.K.,
\lq Fourier Transformations for Difference Analogs of the Harmonic
Oscillator',
in: Proceedings of the XV Workshop on High Energy Physics and Field
Theory,
Protvino, Russia, 6--10 July 1992, to appear.

\item{22.} Askey, R., Atakishiyev, N.M., and Suslov, S.K.,
\lq An Analog of the Fourier Transformations for a $q$-Harmonic
Oscillator',
Preprint No. 5611/1, Kurchatov Institute, Moscow, 1993.

\item{23.} Wilson, J.A., \lq Hypergeometric Series Recurrence
Relations and
Some New Orthogonal Functions', Ph. D. Thesis, University of
Wisconsin,
Madison, Wisc., 1978.

\item{24.} Wilson, J.A., {\it SIAM J. Math. Anal.\/} {\bf 22\/},
1147 (1991).

\item{25.} Rahman, M., {\it Canad. J. Math.\/} {\bf 38}, 605
(1986).

\item{26.} Rahman, M., {\it SIAM J. Math. Anal.\/} {\bf 22\/}, 1430
(1991).

\item{27.} Rahman, M. and Suslov, S.K., \lq Classical Biorthogonal
Rational
Functions', Preprint No. 5614/1, Kurchatov Institute, Moscow 1993;
in { \it Methods of Approximation Theory in Complex Analysis and
Mathematical Physics\/} ( A.A. Gonchar and E.B. Saff, eds. ),
Lecture Notes in Mathematics, Vol. 1550, Springer-Verlag, Berlin,
1993, p. 131.

\item{28.} Andrews, G. and Askey, R., in {\it Polyn\^omes
orthogonaux et
applications\/}, Lecture Notes in Mathematics, Vol. 1171,
Springer-Verlag,
Berlin, 1985, p. 36.

\vfill
\eject

\bye